\documentclass[12pt]{article}%
\usepackage{latexsym}
\usepackage{amsthm}
\usepackage{amsmath}
\usepackage{array}
\usepackage{amsfonts}
\usepackage{amssymb}
\usepackage{graphicx}%
\newtheorem{theorem}{Theorem}
\newtheorem{lemma}{Lemma}

\newtheorem{corollary}{Corollary}
\theoremstyle{definition}

\DeclareMathOperator{\sgn}{sgn}

\begin{document}

\title{A fixed point approach for decaying solutions of difference equations}
\author{Zuzana Do\v{s}l\'{a}, Mauro Marini and Serena Matucci}
\date{ \quad }
\maketitle
\noindent\textbf{Abstract. }{\small A boundary value problem associated to the
difference equation with advanced argument%
\begin{equation}
\label{*}\Delta\bigl (a_{n}\Phi(\Delta x_{n})\bigr)+b_{n}\Phi(x_{n+p}%
)=0,\ \ n\geq1 \tag{$*$}%
\end{equation}
is presented, where $\Phi(u)=|u|^{\alpha}$sgn $u,$ $\alpha>0,p$ is a positive
integer and the sequences $a,b,$ are positive. We deal with a particular type
of decaying solutions of (\ref{*}), that is the so-called intermediate
solutions (see below for the definition) . In particular, we prove the
existence of these type of solutions for (\ref{*}) by reducing it to a
suitable boundary value problem associated to a difference equation without
deviating argument. Our approach is based on a fixed point result for
difference equations, which originates from existing ones stated in the
continuous case. Some examples and suggestions for future researches complete
the paper. }



\vskip4mm

\section{Introduction}

Consider the equation
\begin{equation}
\Delta\bigl (a_{n}\Phi(\Delta x_{n})\bigr)+b_{n}\Phi(x_{n+p})=0,
\tag{P}\label{P}%
\end{equation}
where $\Delta$ is the forward difference operator $\Delta x_{n}=x_{n+1}
-x_{n},$ $\Phi$ is the operator $\Phi(u)=|u|^{\alpha} \sgn u$, $\alpha>0$, $p$
is a positive integer and the sequences $a,b$ are positive for $n\geq1$, and
satisfy
\begin{equation}
\sum_{i=1}^{\infty} \Phi^{\ast}\left(  \frac{1}{a_{i}}\right)  <\infty,
\quad\sum_{i=1}^{\infty} b_{i}=\infty, \label{Hp1}%
\end{equation}
where $\Phi^{\ast}$ is the inverse of the map $\Phi,$ that is $\Phi^{\ast
}(u)=|u|^{1/\alpha}$sgn $u.$

Equation (\ref{P}) appears in the discretization process for searching
spherically symmetric solutions of certain nonlinear elliptic equations with
weighted $\varphi$-Laplacian, see, e.g., \cite{diaz}. A special case of
(\ref{P}) is the discrete half-linear equation
\begin{equation}
\Delta\bigl (a_{n}\Phi(\Delta y_{n})\bigr)+b_{n}\Phi(y_{n+1})=0,
\tag{H}\label{H}%
\end{equation}
which has been studied extensively from various points of view, especially
with regard to the oscillation and the qualitative behavior of nonoscillatory
solutions, see \cite[Chapter 3]{dosly} and references therein.

As usual, a nontrivial solution $x$ of (\ref{P}) is said to be
\textit{nonoscillatory} if $x_{n}$ is either positive or negative for any
large $n$ and \textit{oscillatory} otherwise. In virtue of the Sturm
separation criterion, see, e.g. \cite[Section 8.2.1.]{dosly}, all the
solutions of (\ref{H}) have the same behavior with respect to the oscillation.
In other words, either all nontrivial solutions of (\ref{H}) are
nonoscillatory or all the solutions of (\ref{H}) are oscillatory. Thus,
equivalently, we say that (\ref{H}) is nonoscillatory if (\ref{H}) has a
nonoscillatory solution. Nevertheless, for equation (\ref{P}) with $p\neq1,$
the situation is different, since in this case oscillatory solutions and
nonoscillatory solutions may coexist. \vskip2mm

In this paper we deal with the existence of a particular type of
nonoscillatory solutions, that is solutions $x$ of (\ref{P}) such that
$x_{n}>0,\Delta x_{n}<0$ for large $n$ and
\begin{equation}
\lim_{n}x_{n}=0, \ \ \lim_{n}x_{n}^{[1]}=a_{n}\Phi(\Delta x_{n})=-\infty,
\label{int}%
\end{equation}
where $x^{[1]}$ is called the quasidifference of $x.$ Solutions satisfying
(\ref{int}) are called \textit{intermediate solutions}. This terminology
originates from the corresponding continuous case and it is due to the fact
that, when (\ref{Hp1}) is satisfied, then any nonoscillatory solution $x$ of
(\ref{P}) satisfies either $\lim_{n}x_{n}=\ell_{x}\neq0,$ or (\ref{int}), or
$\lim_{n}x_{n}=0,\lim_{n}x_{n}^{[1]}=-\ell_{x},0<\ell_{x}<\infty,$ see
\cite{novac,kyoto}. The investigation of intermediate solutions is an hard
problem, due to difficulties in finding suitable sharp upper and lower bounds
for these solutions, see, e.g., \cite[page 241]{AgGrOR} and \cite[page 3]%
{KU2}, in which these facts are pointed out for the continuous case.\vskip2mm

In the half-linear case the problem of existence of intermediate solutions has
been completely solved by the following.

\begin{theorem}
\label{Th Half} Assume (\ref{Hp1}).

$(i_{1})$ Equation (\ref{H}) does not have intermediate solutions if
\[
\sum_{n=1}^{\infty}b_{n}\Phi\left(  \sum_{k=n+1}^{\infty}\Phi^{\ast}\left(
\frac{1}{a_{k}}\right)  \right)  +\sum_{n=1}^{\infty}\Phi^{\ast}\left(
\frac{1}{a_{n+1}}\sum_{k=1}^{n}b_{k}\right)  \,<\infty\,.
\]

$(i_{2})$ Equation (\ref{H}) has intermediate solutions if (\ref{H}) is
nonoscillatory and
\[
\sum_{n=1}^{\infty}b_{n}\Phi\left(  \sum_{k=n+1}^{\infty}\Phi^{\ast}\left(
\frac{1}{a_{k}}\right)  \right)  +\sum_{n=1}^{\infty}\Phi^{\ast}\left(
\frac{1}{a_{n+1}}\sum_{k=1}^{n}b_{k}\right)  \,=\infty.
\]

\end{theorem}

\noindent\textit{Proof. } The assertion follows from Theorem 3.1 (b), (c) and
Theorem 3.2. in \cite{summer}, with minor changes. \hfill$\Box$ \vskip4mm

Here, we present a comparison result which allows to solve the boundary value
problem [BVP]
\begin{equation}%
\begin{cases}
\Delta\bigl (a_{n}\Phi(\Delta x_{n})\bigr)+b_{n}\Phi(x_{n+p})=0, & p>1,\\
\lim_{n}x_{n}=0,\ \lim_{n}x_{n}^{[1]}=-\infty &
\end{cases}
\label{BVP}%
\end{equation}
by reducing it to the existence of intermediate solutions in the half-linear
case. The main result is the following.

\begin{theorem}
\label{Th Comp} Assume (\ref{Hp1}) and
\begin{equation}
\limsup_{n} b_{n}<\infty. \label{Hp b}%
\end{equation}
Then the BVP (\ref{BVP}) is solvable, i.e. equation (\ref{P}) with $p>1$ has
intermediate solutions, if and only if the half-linear equation
\begin{equation}
\Delta\bigl (a_{n+p-1}\Phi(\Delta y_{n})\bigr)+b_{n}\Phi(y_{n+1})=0,
\tag{H1}\label{H1}%
\end{equation}
has intermediate solutions.
\end{theorem}

Theorem \ref{Th Comp} will be proved by means of a fixed point result for
discrete operators acting in Fr\'{e}chet spaces, see Theorem \ref{Th fixed} below.

Combining Theorem \ref{Th Half} and the comparison Theorem \ref{Th Comp}, we
get necessary and sufficient conditions for existence of intermediate
solutions of difference equation with advanced argument.

\begin{corollary}
\label{Cor SN} Assume (\ref{Hp1}) and (\ref{Hp b}).

$(i_{1})$ Equation (\ref{P}) with $p>1$ has intermediate solutions if
(\ref{H1}) is nonoscillatory and
\begin{equation}
\sum_{n=1}^{\infty}b_{n}\Phi\left(  \sum_{k=n+1}^{\infty}\Phi^{\ast}\left(
\frac{1}{a_{k+p-1}}\right)  \right)  +\sum_{n=1}^{\infty}\Phi^{\ast}\left(
\frac{1}{a_{n+p}}\sum_{k=1}^{n}b_{k}\right)  \,=\infty\,. \label{S1}%
\end{equation}

(i$_{2})$ Equation (\ref{P}) with $p>1$ does not have intermediate solutions
if
\[
\sum_{n=1}^{\infty}b_{n}\Phi\left(  \sum_{k=n+1}^{\infty}\Phi^{\ast}\left(
\frac{1}{a_{k+p-1}}\right)  \right)  +\sum_{n=1}^{\infty}\Phi^{\ast}\left(
\frac{1}{a_{n+p}}\sum_{k=1}^{n}b_{k}\right)  \,<\infty\,.
\]

\end{corollary}

\vskip6mm

\section{Fixed point approaches}

Boundary value problems for difference equations in $\mathbb{R}^{n}$ are often
solved by reducing the problem to a fixed point equation for a possibly
nonlinear operator in a suitable function space. Thus, the existence of a
solution is obtained by applying a fixed point theorem, for instance\ the
Tychonoff theorem, the Schauder theorem, the Leray-Schauder continuation
principle or Krasnoselkii types fixed point theorems on cones. For a survey on
this topic we refer to the papers \cite{Ag OR,AG,MMR2} and the monographies
\cite{ag-oreg,OR-E}. In particular, in \cite[Chapter 2]{OR-E}, see also
\cite[Chapter 5]{ag-oreg}, certain BVPs are studied by means of a nonlinear
Leray-Schauder alternative. This approach is based on a very general method
given in \cite{FP}. In particular, in \cite{FP} the authors present a
Leray-Schauder continuation principle in locally convex topological vector
spaces, which unifies the Leray-Schauder alternative theorem and the Tichonov
fixed point theorem. More precisely, let $\mathbb{E}$ be a Hausdorff locally
convex topological vector space with a family of seminorms generating the
topology. The following holds.

\begin{theorem}
\label{Th FP} \textrm{\cite[Theorem 1.1]{FP}} 
Let $Q$ be a convex closed
subset of $\mathbb{E}$ and let $T:Q\times\lbrack0,1]\rightarrow\mathbb{E}$ be
a continuous map with relatively compact image. Assume that:

i$_{1})$ $T(x,0)\in Q$ for any $x\in Q;$

i$_{2})$ for any $(x,\lambda)\in\partial Q\times\lbrack0,1)$ with
$T(x,\lambda)=x$ there exists open neighborhoods $U_{x}$ of $x$ in
$\mathbb{E}$ and $I_{\lambda}$ of $\lambda$ in $[0,1)$ such that
\[
T\Bigl((U_{x}\cap\partial Q)\times I_{\lambda}\Bigr)\subset Q.
\]
Then the equation
\[
x=T(x,1)
\]
has a solution.
\end{theorem}

\noindent

Some of the above quoted results have a discrete counterpart. For instance, a
method for solving BVPs associated to difference systems is given in
\cite[Theorem 2.1]{MMR2}. Due to the peculiarities of the discrete case, it
may be applied to functional difference equations, including equations with
deviating arguments or sum difference equations.

\vskip2mm

Now, we present an existence result which generalizes, in the particular case
of scalar difference equations, \cite[Theorem 2.1]{MMR2}.

Denote by $\mathbb{N}_{n}$ and $\mathbb{N}_{m,n}$, the sets
\[
\mathbb{N}_{n}=\{i\in\mathbb{N}:\,i\geq n\in\mathbb{N}\}
\]
\[
\mathbb{N}_{m,n}=\{i\in\mathbb{N}_{m}:\,i<n, \, m,n\in\mathbb{N},m<n\}
\]
and let $\mathbb{X}$ be the space of all real sequences defined on
$\mathbb{N}_{m}$. Hence $\mathbb{X}$ is a Frech\'{e}t space with the topology
of pointwise convergence on $\mathbb{N}_{m}$. From the discrete
Arzel\`{a}-Ascoli theorem (see, e.g., \cite[Theorem 5.3.1]{ag-oreg}), any
bounded set in $\mathbb{X}$ is relatively compact. We recall that a set
$\Omega\subset\mathbb{X}$ is bounded if and only if it consists of sequences
which are equibounded on $\mathbb{N}_{m,n}$ for any $n>m.$ Clearly, if
$\Omega\subset$ $\mathbb{X}$ is bounded, then $\Omega^{\Delta}=\{\Delta
u,\,u\in\Omega\}$ is bounded, too.

\vskip2mm

Using, with minor changes, a discrete counterpart of a compactness and
continuity result stated in \cite[Theorem 1.3]{CFM} for the continuous case,
we have the following.

\begin{theorem}
\label{Th fixed} Consider the BVP
\begin{equation}%
\begin{cases}
\Delta(a_{n}\Phi(\Delta x_{n}))=g(n,x), & n\in\mathbb{N}_{m}\\
x\in S, &
\end{cases}
\label{DF}%
\end{equation}
where $g:\mathbb{N}_{m}\times\mathbb{X}\rightarrow\mathbb{R}$ is a continuous
map, and $S$ is a subset of $\ \mathbb{X}$. \newline Let $G:\,\mathbb{N}
_{m}\times\mathbb{X}^{2}\rightarrow\mathbb{R}$ be a continuous map such that
$G(k,u,u)=g(k,u)$ for all $(k,u)\in\mathbb{N}_{m}\times\mathbb{X}$. If there
exist a nonempty, closed, convex set $\Omega\subset\mathbb{X}$, and a bounded,
closed subset $S_{C}\subset S\cap\Omega$ such that the problem
\begin{equation}%
\begin{cases}
\Delta(a_{n}\Phi(\Delta x_{n}))=G(n,x,q), & n\in\mathbb{N}_{m}\\
x\in S_{C}, &
\end{cases}
\label{DF1}%
\end{equation}
has a unique solution for any $q\in\Omega$ fixed, then (\ref{DF}) has at least
a solution.
\end{theorem}

\noindent\textit{Proof. } Let $T$ be the operator $T:\Omega\rightarrow S_{C}$
which maps every $q\in\Omega$ into the unique solution $x=T(q)$ of
(\ref{DF1}). Let us show that the operator $T$ is continuous with relatively
compact image. The relatively compactness of $T(\Omega)$ follows immediately
since $S_{C}$ is bounded. To prove the continuity of $T$ in $\Omega$, let
$\{q^{j}\}$ be a sequence in $\Omega$, ${q}^{j}\rightarrow{q}^{\infty}
\in\Omega$, and let ${v}^{j}=T(q^{j})$. Since $T(\Omega)$ is relatively
compact, $\{{v}^{j}\}$ admits a subsequence (still indicated with $\{{v}
^{j}\}$) which converges to ${v}^{\infty}\in\mathbb{X}$. As $v^{j}\in S_{C}$
and $S_{C}$ is closed, then ${v}^{\infty}\in S_{C}$. Taking into account the
continuity of $G$, we obtain
\[
\Delta(a_{n}\Phi(\Delta v_{k}^{\infty}))=\lim_{j}\Delta(a_{n}\Phi(\Delta
v_{k}^{j}))=\lim_{j}G(k,q^{j},v^{j})=G(k,q^{\infty},{v}^{\infty}).
\]
The uniqueness of the solution of (\ref{DF1}) yields $v^{\infty}=T({q}
^{\infty})$, and therefore $T$ is continuous on $\Omega$. By the Tychonoff
fixed point theorem, $T$ has at least one fixed point in $\Omega$, which is a
solution of (\ref{DF}), as it can be easily checked, taking into account that
$S_{C}\subset S$. \hfill$\Box$ \vskip4mm

As follows from the proof of Theorem \ref{Th fixed}, no explicit form of the
fixed point operator is needed for the solvability of (\ref{DF}). A key point
for the unique solvability of (\ref{DF1}) in Theorem \ref{Th fixed} is the
choice of the map $G$. To this aim, in our opinion, the best cases are the
following two, namely
\[
i_{1})\text{ \ }G(n,q,x)=\tilde{g}(n,q),\text{ \ \ \ }i_{2})\text{
\ \ }G(n,q,x)=\tilde{g}(n,q)\Phi(x_{n+1}).
\]
In the case $i_{2}),$ the equation in (\ref{DF1}) is a half-linear equation
and in this situation a very large variety of results is known, see, e.g.,
\cite[Chapter VIII]{dosly}. An application in this direction is in
\cite[Section 4]{MMR2}.

In the case $i_{1}),$ \ that is when the function $G$ does not depend on $x$
and the equation in (\ref{DF1}) is affine, Theorem \ref{Th fixed} can be
particularly useful to solve BVPs associated to difference equations with
deviating arguments. Indeed, in this case, it can lead to a BVP associated to
a second order difference equation without deviating argument. An application
of this fact is in the following section. \vskip2mm

\section{Proof of Theorem \ref{Th Comp}}

For proving Theorem \ref{Th Comp}, the following auxiliary result is needed.

\begin{lemma}
\label{Lem Est} Assume (\ref{Hp1}) and (\ref{Hp b}). Let $x$ be an eventually
positive decreasing solution of (\ref{P}) such that $\lim_{n}x_{n}=0$. Then
the series
\[
\sum_{i=2}^{\infty}\Phi^{\ast}\left(  \frac{1}{a_{i+p-1}}\sum_{k=1}^{i-1}
b_{k}\Phi(x_{k+p})\right)
\]
converges.
\end{lemma}

\noindent\textit{Proof. } Without loss of generality, suppose $0<x_{n}\leq1,$
$\Delta x_{n}<0$ for $n\geq m_{0}\geq1$. We claim that for any $k,j\geq m_{0}$
we have
\begin{equation}
\left\vert x_{j}^{[1]}-x_{k}^{[1]}\right\vert \leq B|j-k|, \label{a:20}%
\end{equation}
where $B=\sup_{i\geq1}b_{i}.$ For simplicity, suppose $k\geq j$. Summing
(\ref{P}) we have
\begin{equation}
\label{x2}x_{k}^{[1]}-x_{j}^{[1]}=-\sum_{i=j}^{k-1}b_{i}\Phi(x_{i+p}),
\end{equation}
where, as usual, $\sum_{i=k_{1}}^{k}\gamma_{i}=0$ if $k<k_{1}$. Since
$\left\vert x_{i}\right\vert \leq1$ and $b_{i}\leq B,$ the inequality
(\ref{a:20}) follows. Further, from (\ref{x2}),
\begin{gather*}
\sum_{i=m_{0}}^{\infty}\Phi^{\ast}\left(  \frac{1}{a_{i+p-1}}\sum_{k=m_{0}%
}^{i-1} b_{k}\Phi(x_{k+p})\right)  = \sum_{i=m_{0}}^{\infty}\Phi^{\ast}\left(
\frac{1}{a_{i+p-1}}\left(  x_{m_{0}}^{[1]}-x_{i}^{[1]}\right)  \right)  =\\
=\sum_{i=m_{0}}^{\infty}\Phi^{\ast}\left(  \frac{1}{a_{i+p-1}}\left(
x_{m_{0}}^{[1]}- x_{i}^{[1]}-x_{i+p-1}^{[1]}+x_{i+p-1}^{[1]}\right)  \right)
.
\end{gather*}
Since
\[
\left\vert x_{i}^{[1]}-x_{i+p-1}^{[1]}+x_{i+p-1}^{[1]}-x_{m_{0}}
^{[1]}\right\vert \leq\left\vert x_{i}^{[1]}-x_{i+p-1}^{[1]}\right\vert
+\left\vert x_{i+p-1}^{[1]}\right\vert +\left\vert x_{m_{0}}^{[1]}\right\vert
,
\]
in view of (\ref{a:20}) we obtain
\begin{gather*}
\left\vert \sum_{i=m_{0}}^{\infty}\Phi^{\ast}\left(  \frac{1}{a_{i+p-1}
}\left(  x_{i}^{[1]}-x_{m_{0}}^{[1]}\right)  \right)  \right\vert \leq\\
\leq\sum_{i=m_{0}}^{\infty}\Phi^{\ast}\left(  \frac{B(p-1)+\left\vert
x_{m_{0}}^{[1]}\right\vert }{a_{i+p-1}}+\frac{\left\vert x_{i+p-1}
^{[1]}\right\vert }{a_{i+p-1}}\right)  .
\end{gather*}
In virtue of (\ref{Hp1}), the series
\[
\sum_{i=m_{0}}^{\infty}\Phi^{\ast}\left(  \frac{B(p-1)+\left\vert x_{m_{0}
}^{[1]}\right\vert }{a_{i+p-1}}\right)
\]
converges. Since
\[
\sum_{i=m_{0}}^{\infty}\Phi^{\ast}\left(  \frac{\left\vert x_{i+p-1}
^{[1]}\right\vert }{a_{i+p-1}}\right)  =-\sum_{i=m_{0}}^{\infty}\Delta
x_{i+p-1}=x_{m_{0}+p-1},
\]
using (\ref{Hp1}) and the inequality
\[
\Phi^{\ast}(X+Y)\leq\sigma_{\alpha}(\Phi^{\ast}(X)+\Phi^{\ast}(Y)),
\]
where
\[
\sigma_{\alpha}=\left\{
\begin{array}
[c]{ccc}%
1 &  & \text{if }\alpha\geq1\\
2^{(1-\alpha)/\alpha} &  & \text{if }\alpha<1
\end{array}
\right.  ,
\]
we obtain the assertion. \hfill$\Box$ \vskip4mm

\textit{Proof of Theorem \ref{Th Comp}. }First, we prove that if (\ref{P}) has
intermediate solutions, then (\ref{H1}) has intermediate solutions.

\vskip2mm

Let $x$ be an intermediate solution of (\ref{P}) and, without loss of
generality, assume for $n\geq n_{0}\geq1$
\begin{equation}
0<x_{n}<1,\text{ \ \ }\Delta x_{n}<0\text{ }. \label{a:10}%
\end{equation}
In view of (\ref{Hp b}), there exists $L>0$ such that for any $n\geq n_{0}$
\begin{equation}
\sum_{i=n}^{n+p-2}b_{i}\leq L. \label{a:14}%
\end{equation}
Moreover, let $M$ be a positive constant, $M<1,$ such that
\begin{equation}
\Phi(M)\leq\frac{\left\vert x_{n_{0}}^{[1]}\right\vert }{L+\left\vert
x_{n_{0}}^{[1]}\right\vert }. \label{a:16}%
\end{equation}

Let $\mathbb{X}$ be the Fr\'{e}chet space of real sequences defined for $n\geq
n_{0}$, endowed with the topology of convergence on $\mathbb{N}_{n_{0} }$, and
consider the subset $\Omega\subset$ $\mathbb{X}$ defined by
\[
\Omega=\big\{u\in\mathbb{X}\,:Mx_{n+p-1}\leq u_{n}\leq x_{n+p-1}\big\}.
\]
For any $u\in\Omega$ consider the boundary value problem [BVP]
\begin{equation}%
\begin{cases}
\Delta\bigl (a_{n+p-1}\Phi(\Delta z_{n})\bigr)+b_{n}\Phi(u_{n+1})=0, & n \geq
n_{0}\\
z_{n_{0}}^{[1]}=x_{n_{0}}^{[1]}, \ \lim_{n}z_{n}=0, &
\end{cases}
\label{BVP lin}%
\end{equation}
where $z^{[1]}$ denotes the quasidifference of $z,$ that is
\begin{equation}
z_{n}^{[1]}=a_{n+p-1}\Phi(\Delta z_{n}). \label{z1}%
\end{equation}
For any $u\in\Omega$ we have
\[
\sum_{k=n_{0}}^{n}b_{k}\Phi(u_{k+1})\leq\sum_{k=n_{0}}^{n}b_{k}\Phi(x_{k+p}).
\]
Hence, using Lemma \ref{Lem Est}, we have
\[
\lim_{n}\sum_{i=n}^{\infty}\Phi^{\ast}\left(  \frac{1}{a_{i+p-1}}\left(
\sum_{k=n_{0}}^{i-1}b_{k}\Phi(u_{k+1})\right)  \right)  =0.
\]
Thus, a standard calculation shows that for any $u\in\Omega$ the BVP
(\ref{BVP lin}) has the unique solution $z$. Let $\mathcal{T}$ \ be the
operator which associates to any $u\in\Omega$ the unique solution $z$ of
(\ref{BVP lin}).

\medskip\noindent Summing the equation in (\ref{BVP lin}) and using (\ref{x2})
we get
\[
z_{n}^{[1]}=x_{n_{0}}^{[1]}-\sum_{k=n_{0}}^{n-1}b_{k}\Phi(u_{k+1})\geq
x_{n_{0}}^{[1]}-\sum_{k=n_{0}}^{n-1}b_{k}\Phi(x_{k+p})=x_{n}^{[1]}.
\]
Since $x^{[1]}$ is decreasing for $n\geq n_{0}$ and $p>1,$ we obtain for
$n\geq n_{0}$
\begin{equation}
z_{n}^{[1]}\geq x_{n+p-1}^{[1]} \label{a:80}%
\end{equation}
i.e., in view of (\ref{z1}),
\[
a_{n+p-1}\Phi(\Delta z_{n})\geq a_{n+p-1}\Phi(\Delta x_{n+p-1}),
\]
that is,
\[
\Delta z_{n}\geq\Delta x_{n+p-1}.
\]
Since $\lim_{i}z_{i}=\lim_{i}x_{i}=0,$ we obtain for $n\geq n_{0}$
\begin{equation}
z_{n}\leq x_{n+p-1}. \label{a:100}%
\end{equation}
Now, let us prove that for\ $n\geq n_{0}$
\begin{equation}
z_{n}\geq Mx_{n+p-1}. \label{a:110}%
\end{equation}
Summing the equation in (\ref{BVP lin}) and using (\ref{z1}) we get
\[
z_{n}^{[1]}=x_{n_{0}}^{[1]}-\sum_{k=n_{0}}^{n-1}b_{k}\Phi(u_{k+1})\leq
x_{n_{0}}^{[1]}-\Phi(M)\sum_{k=n_{0}}^{n-1}b_{k}\Phi(x_{k+p}),
\]
or, using (\ref{x2}),
\begin{align}
z_{n}^{[1]}  &  \leq x_{n_{0}}^{[1]}+\Phi(M)\left(  x_{n}^{[1]}-x_{n_{0}%
}^{[1]}\right)  =\nonumber\\
&  \text{ }\label{a:41}\\
&  =\Phi(M)x_{n+p-1}^{[1]}+\Phi(M)\left(  x_{n}^{[1]}-x_{n+p-1}^{[1]}\right)
+(1-\Phi(M))x_{n_{0}}^{[1]}.\nonumber
\end{align}
From (\ref{x2}), (\ref{a:10}) and (\ref{a:14}), we also have
\[
x_{n}^{[1]}-x_{n+p-1}^{[1]}=\sum_{i=n}^{n+p-2}b_{i}\Phi(x_{i+p})\leq\sum
_{i=n}^{n+p-2}b_{i}\leq L.
\]
Thus, from (\ref{a:41}) we obtain
\begin{equation}
z_{n}^{[1]}\leq\text{ }\Phi(M)x_{n+p-1}^{[1]}+\left(  L+\left\vert x_{n_{0}%
}^{[1]}\right\vert \right)  \Phi(M)+x_{n_{0}}^{[1]}. \label{a:63}%
\end{equation}
In view of (\ref{a:16}) we have
\[
\left(  L+\left\vert x_{n_{0}}^{[1]}\right\vert \right)  \Phi(M)+x_{n_{0}%
}^{[1]}\leq0.
\]
Hence, from (\ref{a:63}) we get for $n\geq n_{0}$
\begin{equation}
z_{n}^{[1]}\leq\Phi(M)x_{n+p-1}^{[1]} \label{a:90}%
\end{equation}
or, in view of (\ref{z1}), $\Delta z_{n}\leq M\text{ }\Delta x_{n+p-1},$ and
(\ref{a:110}) follows, since $\lim_{i}z_{i}=\lim_{i}x_{i}=0.$ Thus, in virtue
of (\ref{a:100}) and (\ref{a:110}), the operator $\mathcal{T}$ maps $\Omega$
into itself, that is%
\[
\mathcal{T}(\Omega)\subset\Omega.
\]

Denote by $S$ the boundary conditions in (\ref{BVP lin}), i.e.
\[
S=\left\{  v\in\mathbb{X}:\,a_{n_{0}+p-1}\Phi(\Delta v_{n_{0}})=x_{n_{0}%
}^{[1]},\ \lim_{n}v_{n}=0\right\}
\]
For any $z\in\mathcal{T}(\Omega)$ we have $z\in S$. Since $\mathcal{T}%
(\Omega)\subset\Omega$, we get $z\in\Omega\cap S$.

Denote by $S_{C}$ the subset of $\mathbb{X}$ given by
\[
S_{C}=S\cap\Omega.
\]
Since $\lim_{n}x_{n}=0$, it holds
\[
S_{C}=\left\{  v\in\mathbb{X}:\,a_{n_{0}+p-1}\Phi(\Delta v_{n_{0}})=x_{n_{0}%
}^{[1]},\,Mx_{n+p-1}\leq v_{n}\leq x_{n+p-1}\right\} \, .
\]
Thus $S_{C}$ is a bounded and closed subset of $\mathbb{X.}$ Applying Theorem
\ref{Th fixed} we obtain that the operator $\mathcal{T}$ \ has a fixed point
$\overline{z}\in S_{C}$. \ Clearly the sequence $\overline{z}$ is a solution
of (\ref{H1}) and $\lim_{n}\overline{z}_{n}=0.$ Since $\overline{z}%
\in\mathcal{T}(\Omega),$ from (\ref{a:80}) and (\ref{a:90}), we get
\[
x_{n+p-1}^{[1]}\leq\overline{z}_{n}^{[1]}\leq\Phi(M)x_{n+p-1}^{[1]}%
\]
and so $\lim_{n}\overline{z}_{n}^{[1]}=-\infty.$ Hence $\overline{z}$ is an
intermediate solution of (\ref{H1}).

\vskip2mm

Now, we prove the vice-versa, that is if\ (\ref{H1}) has intermediate
solutions, then (\ref{P}) has intermediate solutions.

\vskip2mm

The argument is similar to the one above given, with minor changes. Let $y$ be
an intermediate solution of (\ref{H1}) such that for $n\geq n_{0}\geq1$
\[
0<y_{n}<1,\text{ \ \ }\Delta y_{n}<0\text{,}%
\]
and define
\[
n_{1}=n_{0}+p.
\]
In view of (\ref{Hp b}), there exists $\Lambda>0$ such that for any $n\geq
n_{1}$
\begin{equation}
\sum_{i=n-p+1}^{n-1}b_{i}\leq\Lambda. \label{a:170}%
\end{equation}
Without loss of generality, we can suppose
\begin{equation}
\Lambda<\left\vert y_{n_{1}}^{[1]}\right\vert , \label{a:172}%
\end{equation}
where $y_{n_{1}}^{[1]}=a_{n_{1}+p-1}\Phi(\Delta y_{n_{1}}).$ Moreover, let
$H>1$ be a positive constant such that
\begin{equation}
\Phi(H)\geq\frac{\left\vert y_{n_{1}}^{[1]}\right\vert }{\left\vert y_{n_{1}%
}^{[1]}\right\vert -\Lambda}. \label{a:173}%
\end{equation}
Let $\mathbb{X}_{1}$ be the Fr\'{e}chet space of the real sequences defined
for $n\geq n_{1},$ endowed with the topology of convergence on $\mathbb{N}%
_{n_{1}}$. Define the subset $\Omega_{1}$ of $\mathbb{X}_{1}$
\[
\Omega_{1}=\big\{u\in\mathbb{X}_{1}\,:y_{n-p+1}\leq u_{n}\leq Hy_{n-p+1}%
\big\}
\]
and for any $u\in\Omega_{1}$ consider the BVP
\begin{equation}%
\begin{cases}
\Delta\bigl (a_{n}\Phi(\Delta w_{n})\bigr)+b_{n}\Phi(u_{n+p})=0, & n\geq
n_{1}\\
w_{n_{1}}^{[1]}=y_{n_{1}}^{[1]},\ \lim_{n}w_{n}=0, &
\end{cases}
\label{BVP lin2}%
\end{equation}
where $w^{[1]}$ and $y^{[1]}$ denote the quasidifferences of $w$ and $y,$
respectively, that is the sequences
\begin{equation}
w_{n}^{[1]}=a_{n}\Phi(\Delta w_{n}),\text{ \ \ \ }y_{n}^{[1]}=a_{n+p-1}%
\Phi(\Delta y_{n}).\text{\ } \label{N:10}%
\end{equation}
As before, for any $u\in\Omega_{1}$ the BVP (\ref{BVP lin2}) has a unique
solution $w$. Thus, let $\mathcal{T}$ be the operator which associates to any
$u\in\Omega_{1}$ the unique solution $w$ of (\ref{BVP lin2}). Summing the
equation in (\ref{BVP lin2}) and using (\ref{N:10}) we get
\begin{equation}
w_{n}^{[1]}=y_{n_{1}}^{[1]}-\sum_{k=n_{1}}^{n-1}b_{k}\Phi(u_{k+p})\leq
y_{n_{1}}^{[1]}-\sum_{k=n_{1}}^{n-1}b_{k}\Phi(y_{k+1})=y_{n}^{[1]}\,.
\label{N:15}%
\end{equation}
Since $y^{[1]}$ is decreasing for $n\geq n_{0}$ and $p>1,$ we have
$y_{n}^{[1]}\leq y_{n-p+1}^{[1]}$ for $n\geq n_{1}$, and from (\ref{N:15}) we
obtain
\begin{equation}
w_{n}^{[1]}\leq y_{n-p+1}^{[1]} \label{N:18}%
\end{equation}
i.e.%
\[
\Delta w_{n}\leq\Delta y_{n-p+1},
\]
which implies
\begin{equation}
w_{n}\geq y_{n-p+1}, \label{a:200}%
\end{equation}
since $\lim_{i}w_{i}=\lim_{i}y_{i}=0$.

Now, let us prove that
\begin{equation}
w_{n}\leq H\text{ }y_{n-p+1}\text{ \ \ }. \label{a:250}%
\end{equation}
Summing the equation in (\ref{BVP lin2}) we get
\[
w_{n}^{[1]}=y_{n_{1}}^{[1]}-\sum_{k=n_{1}}^{n-1}b_{k}\Phi(u_{k+p})\geq
y_{n_{1}}^{[1]}-\Phi(H)\sum_{k=n_{1}}^{n-1}b_{k}\Phi(y_{k+1}),
\]
that is
\begin{equation}%
\begin{split}
w_{n}^{[1]}  &  \geq y_{n_{1}}^{[1]}+\Phi(H)\left(  y_{n}^{[1]}-y_{n_{1}%
}^{[1]}\right) \\
&  =\Phi(H)y_{n-p+1}^{[1]}+\left(  y_{n}^{[1]}-y_{n-p+1}^{[1]}\right)
\Phi(H)+(1-\Phi(H))y_{n_{1}}^{[1]}.
\end{split}
\label{a:210}%
\end{equation}
From (\ref{H1}) and (\ref{a:170}) we have
\[
y_{n}^{[1]}-y_{n-p+1}^{[1]}=-\sum_{i=n-p+1}^{n-1}b_{i}\Phi(y_{i+1}%
)\geq-\Lambda.
\]
Thus, from (\ref{a:210}) we obtain
\begin{align*}
w_{n}^{[1]}  &  \geq\Phi(H)y_{n-p+1}^{[1]}-\Lambda\Phi(H)+y_{n_{1}}^{[1]}%
-\Phi(H)y_{n1}^{[1]}\\
&  =\Phi(H)y_{n-p+1}^{[1]}+\left(  \left\vert y_{n_{1}}^{[1]}\right\vert
-\Lambda\right)  \Phi(H)+y_{n_{1}}^{[1]},
\end{align*}
or, in view of (\ref{a:172}) and (\ref{a:173}),
\begin{equation}
w_{n}^{[1]}\geq\Phi(H)y_{n-p+1}^{[1]} \label{a:230}%
\end{equation}
i.e.,
\[
\Delta w_{n}\geq H\text{ }\Delta y_{n-p+1}.
\]
Since $\lim_{i}w_{i}=\lim_{i}y_{i}=0,$ from here we get (\ref{a:250}). Hence,
in virtue of (\ref{a:200}) and (\ref{a:250}), the operator $\mathcal{T}$ maps
$\Omega_{1}$ into itself, i.e.%
\[
\mathcal{T}(\Omega_{1})\subset\Omega_{1}.
\]

\medskip Using the same argument to the one given in the sufficiency part,
denote by $S_{1}$ the boundary conditions in (\ref{BVP lin2}). Applying
Theorem \ref{Th fixed} with $S_{C}=S_{1}\cap\Omega_{1}$, we get that the
operator $\mathcal{T}$ \ has a fixed point $\overline{w}\in S_{1}\cap
\Omega_{1}.$ Clearly the sequence $\overline{w}$ is a solution of (\ref{P})
and $\lim_{n}\overline{w}_{n}=0.$ Since $\overline{w}\in\mathcal{T}(\Omega
_{1})$, from (\ref{N:18}) and (\ref{a:230}), we get
\[
Hy_{n-p+1}^{[1]}\leq\overline{w}_{n}^{[1]}\leq y_{n-p+1}^{[1]}%
\]
and so $\lim_{n}\overline{w}_{n}^{[1]}=-\infty.$

Hence $\overline{w}$ is an intermediate solution of (\ref{P}) and the proof is
complete.\hfill$\Box$ \vskip2mm

\section{Suggestions and examples}

The following example illustrates Theorem \ref{Th Comp} and Corollary
\ref{Cor SN}.\vskip2mm

\noindent\textbf{Example 1. }Consider the difference equation with advanced
argument
\begin{equation}
\Delta\bigl ((n-p+1)^{1+\alpha}\Phi(\Delta x_{n})\bigr)+\gamma\Phi
(x_{n+p})=0,\text{ \ \ }n\geq p\geq2, \label{Ex1}%
\end{equation}
where $\gamma$ is a positive constant. Using Theorem \ref{Th Comp} and
Corollary \ref{Cor SN}, it is easy to show that (\ref{Ex1}) has intermediate
solutions if and only if
\begin{equation}
0<\gamma\leq\left(  \frac{1}{1+\alpha}\right)  ^{\alpha+1}. \label{crit}%
\end{equation}
Indeed, consider the half-linear equation
\begin{equation}
\Delta\bigl (n^{1+\alpha}\Phi(\Delta x_{n})\bigr)+\gamma\Phi(x_{n+1})=0.
\label{Eul1}%
\end{equation}
A standard calculation shows that (\ref{Hp1}) is satisfied. Moreover, using
the change of variable
\begin{equation}
y_{n}=n^{1+\alpha}\Phi(\Delta x_{n}) \label{qd}%
\end{equation}
the equation (\ref{Eul1}) is transformed into the generalized discrete Euler
equation
\[
\Delta\bigl (\Phi^{\ast}(\Delta y_{n})\bigr)+\gamma^{1/\alpha}\left(  \frac
{1}{n+1}\right)  ^{(1+\alpha)/\alpha}\Phi^{\ast}(y_{n+1})=0,
\]
which is nonoscillatory if $\gamma$ satisfies (\ref{crit}) and oscillatory if
\begin{equation}
\gamma>\left(  \frac{1}{1+\alpha}\right)  ^{\alpha+1}, \label{crit2}%
\end{equation}
see, e.g., \cite{Re1}. Since the transformation (\ref{qd}) maintains the
oscillatory behavior, see, e.g., \cite{d-r2}, we get that (\ref{Eul1}) is
nonoscillatory if and only if (\ref{crit}) is satisfied. Moreover, as
$a_{n}=n^{1+\alpha},$ $b_{n}=\gamma,$ we have
\[
\sum_{n=p}^{\infty}\Phi^{\ast}\left(  \frac{1}{a_{n+p}}\sum_{k=p}^{n}
b_{k}\right)  =\gamma^{1/\alpha}\sum_{n=p}^{\infty}\left(  \frac{n}
{n+p}\right)  ^{1/\alpha}\frac{1}{n+p}=\infty
\]
and so the condition (\ref{S1}) is satisfied. Hence, from Corollary
\ref{Cor SN}$-i_{1})$ equation (\ref{Ex1}) has intermediate solutions if
(\ref{crit}) is satisfied. When (\ref{crit2}) holds, as claimed, the
half-linear equation (\ref{Eul1}) is oscillatory and it does not admit
intermediate solutions. Thus, from Theorem \ref{Th Comp}, equation (\ref{Ex1})
does not have intermediate solutions, too.

\vskip6mm

The existence of intermediate solutions for (\ref{P}) does not depend on
condition (\ref{Hp b}), as the following example shows.\vskip2mm

\noindent\textbf{Example 2. }Consider the difference equation with advanced
argument
\begin{equation}
\Delta\bigl ((n+1)!\text{ }\Delta x_{n}\bigr)+(n+p)!\text{ }x_{n+p}=0,\text{
\ }p\geq2. \label{Ex0}%
\end{equation}
A direct computation shows that $x=\left\{  1/n!\right\}  $ is an intermediate
solution of (\ref{Ex0}). Nevertheless, for (\ref{Ex0}), assumption
(\ref{Hp b}) is not verified and Theorem \ref{Th Comp} cannot be applied.
Hence, it is an open problem if Theorem \ref{Th Comp} continues to hold when
condition (\ref{Hp b}) failed.

\vskip6mm

Example 1 suggests the following two comparison results.

\begin{corollary}
\label{Cor Eul 1} Assume (\ref{Hp1}) and (\ref{Hp b}). Suppose that
\[
a_{n}\geq n^{1+\alpha},\text{ \ \ \ and \ \ }
{\displaystyle\sum\limits_{n=N}^{\infty}}
\Phi^{\ast}\left(  \frac{n}{a_{n+p}}\right)  =\infty,
\]
where $N\geq p.$ Set $L=\sup_{n\geq N}b_{n}.$ If $b$ is bounded away from zero
and
\[
L<\left(  \frac{1}{1+\alpha}\right)  ^{1+\alpha},
\]
then the equation
\begin{equation}
\Delta\bigl ((n-p+1)^{1+\alpha}\Phi(\Delta x_{n})\bigr)+b_{n}\Phi
(x_{n+p})=0,\text{ \ \ }n\geq p\geq2, \label{P2}%
\end{equation}
has intermediate solutions.
\end{corollary}

\noindent\textit{Proof. }Consider the equation
\begin{equation}
\Delta\bigl (n^{1+\alpha}\Phi(\Delta x_{n})\bigr)+L\Phi(x_{n+1})=0. \label{eC}%
\end{equation}
Reasoning as in Example 1, we get that (\ref{eC}) is nonoscillatory. Hence, in
virtue of the Sturm comparison theorem, see, e.g., \cite[Chapter 8.2]{dosly},
the half-linear equation (\ref{H}) is nonoscillatory. Moreover, since $b$ is
bounded away from zero, there exists $\varepsilon>0$ such that $b_{n}
\geq\varepsilon$ for any $n\geq1.$ Hence
\[%
{\displaystyle\sum\limits_{n=N}^{\infty}}
\Phi^{\ast}\left(  \frac{1}{a_{n+1}}\sum_{k=1}^{n}b_{k}\right)  \geq\Phi
^{\ast}(\varepsilon)%
{\displaystyle\sum\limits_{n=N}^{\infty}}
\Phi^{\ast}\left(  \frac{n}{a_{n+1}}\right)  =\infty,
\]
which implies (\ref{S1}). Hence, applying Corollary \ref{Cor SN}$-i_{1})$ to
the equation (\ref{P2}), we get the assertion. \hfill$\Box$ \vskip4mm

\begin{corollary}
\label{Cor Eul 2} Assume (\ref{Hp1}) and (\ref{Hp b}). Moreover, suppose that
for $n\geq N\geq p,$
\[
a_{n}\leq n^{1+\alpha}.
\]
Set $\ell=\inf_{n\geq N}b_{n}.$ If
\[
\ell>\left(  \frac{1}{1+\alpha}\right)  ^{1+\alpha},
\]
then the equation (\ref{P2}) does not have intermediate solutions.
\end{corollary}

\noindent\textit{Proof. }The argument is similar to the one given in Corollary
\ref{Cor Eul 1}. Consider the half-linear equation
\begin{equation}
\Delta\bigl (n^{1+\alpha}\Phi(\Delta x_{n})\bigr)+\ell\text{ }\Phi(x_{n+1})=0.
\label{8b}%
\end{equation}
Reasoning as before, we get that (\ref{8b}) is oscillatory. Hence, in virtue
of the Sturm comparison theorem, also (\ref{H}) is oscillatory, and so
(\ref{H}) does not have intermediate solutions. Thus, applying Theorem
\ref{Th Comp} to the equation (\ref{H}), we obtain the assertion. \hfill$\Box$ \vskip6mm

Some suggestions for future researches are in order.\vskip2mm

\textbf{1)} As claimed, Theorem \ref{Th FP} from \cite{FP} gives a very
general fixed point result, which is based on a continuation principle in a
Hausdorff locally convex space. Further, in \cite{FP} the solvability to
certain BVPs in the continuous case is also given. It should be interesting to
establish corresponding discrete versions of these existence results,
especially for \cite[Theorem 2.2]{FP}, which deals with a scalar equation. It
could be useful for studying discrete BVPs when it is hard to find an
appropriate bounded closed set $\Omega$ which is mapped into itself, as, for
instance occurs for intermediate solutions to Emden-Fowler superlinear
discrete equation
\[
\Delta(a_{n}|\Delta x_{n}|^{\alpha}\operatorname{sgn}\Delta x_{n}%
)+b_{n}|x_{n+1}|^{\beta}\operatorname{sgn}x_{n+1}=0\,,\text{ \ \ }\alpha
<\beta,
\]
see, e.g., \cite[Section V.]{summer}$.$\vskip2mm

2) The proof of Theorem \ref{Th Comp} does not work if $p\leq0.$ Indeed, in
this case the half-linear equation (\ref{H1}) is not defined, due to the shift
in the weight coefficient $a$ of the discrete operator%
\[
(Dx)_{n}=\Delta(a_{n}\Phi(\Delta x_{n})).
\]
Further, when $p<0,$ the sequence $\left\{  u_{n+p}\right\}  $ in
BVP\ (\ref{BVP lin2}) has to be defined not only for $n\geq n_{1},$ but also
for any $i\geq n_{1}+p.$ Consequently, when (\ref{P}) is an equation with
delay, the solvability of (\ref{BVP}) requires a different approach, which
will be presented in a forthcoming paper \cite{DMM}.\vskip2mm

3) Theorem \ref{Th Comp} establishes\ a comparison between the asymptotic
decay of intermediate solutions of (\ref{P}) with the one of an associated
half-linear equation. Recently, in some particular cases, a precise asymptotic
analysis of intermediate solutions for discrete half-linear equations has been
made in the framework of regular variation, see \cite{MR}. It should be
interesting to apply this approach for obtaining a precise description of the
asymptotic behavior also for intermediate solutions of the equations with
advanced argument. 

\vskip6mm

\noindent{\small \emph{Authors' addresses}: }

{\small \noindent\emph{Zuzana Do\v{s}l\'{a},} Department of Mathematics and
Statistics, Masaryk University, Kotl\'{a}\v{r}sk\'{a} 2, CZ-61137 Brno, Czech
Republic. \newline E-mail: \texttt{dosla@math.muni.cz}} \vskip2mm

{\small \noindent\emph{Mauro Marini,} Department of Mathematics and Computer
Science "Ulisse Dini", University of Florence, Via di S. Marta 3, I-50139
Florence, Italy. \newline E-mail: \texttt{mauro.marini@unifi.it}} \vskip2mm

{\small \noindent\emph{Serena Matucci,} Department of Mathematics and Computer
Science "Ulisse Dini", University of Florence, Via di S. Marta 3, I-50139
Florence, Italy. \newline E-mail: \texttt{serena.matucci@unifi.it}}

\end{document}